# ASYMPTOTIC BEHAVIOR OF THE POISSON–DIRICHLET DISTRIBUTION FOR LARGE MUTATION RATE[1]

BY DONALD A. DAWSON AND SHUI FENG

*Carleton University and McMaster University*

The large deviation principle is established for the Poisson–Dirichlet distribution when the parameter $\theta$ approaches infinity. The result is then used to study the asymptotic behavior of the homozygosity and the Poisson–Dirichlet distribution with selection. A phase transition occurs depending on the growth rate of the selection intensity. If the selection intensity grows sublinearly in $\theta$, then the large deviation rate function is the same as the neutral model; if the selection intensity grows at a linear or greater rate in $\theta$, then the large deviation rate function includes an additional term coming from selection. The application of these results to the heterozygote advantage model provides an alternate proof of one of Gillespie's conjectures in [*Theoret. Popul. Biol.* **55** 145–156].

**1. Introduction.** Let

$$\nabla = \left\{ (p_1, p_2, \ldots) : p_1 \geq p_2 \geq \cdots \geq 0, \sum_{k=1}^{\infty} p_k = 1 \right\}.$$

The Poisson–Dirichlet distribution with parameter $\theta > 0$ [henceforth denoted by $PD(\theta)$] is a probability measure on $\nabla$. It was introduced by Kingman [10] as an asymptotic distribution of the order statistics of a symmetric Dirichlet distribution with parameters $K, \alpha$ when $K \to \infty$ and $\alpha \to 0$ in a way such that $\lim_{K \to \infty} K\alpha = \theta$. The distribution coincides with the distribution of the normalized jump sizes of a Gamma process over the interval $(0, \theta)$ ranked in descending order. We use $\mathbf{P}(\theta) = (P_1(\theta), P_2(\theta), \ldots)$ to denote the $\nabla$-valued random variable with distribution $PD(\theta)$. $PD(\theta)$ appears in many different contexts, including Bayesian statistics, number theory, combinatorics and population genetics. In the context of population genetics,

Received April 2004; revised May 2005.
[1]Supported by NSERC.

the distribution describes the equilibrium proportions of different alleles in the infinitely many neutral alleles model. The component $P_k(\theta)$ represents the proportion of the $k$th most frequent allele. If $u$ is the individual mutation rate and $N$ is the effective population size, then the parameter $\theta = 4Nu$ is the population mutation rate. When $\theta$ is small, a large proportion of the population tends to concentrate on a small set of alleles, whereas for large $\theta$, the population is fairly evenly spread. A more friendly way of describing the distribution $PD(\theta)$ is through the following size-biased sampling process. We first let $U_k, k = 1, 2, \ldots$, be a sequence of independent, identically distributed random variables with common distribution $Beta(1, \theta)$. We then generate a random probability vector representing allelic frequencies as follows:

$$X_1 = U_1, \qquad X_n = (1 - U_1) \cdots (1 - U_{n-1}) U_n, \qquad n \geq 2.$$

In other words, the frequency of the first allele type is chosen at random, this is removed and the relative frequency of the second allele is chosen in the same way. This pattern is repeated to get all samples. Then the frequency of allelic type of the $k$th selected sample will be $X_k$. It can be shown that $X_1, X_2, \ldots$ reordered in descending order has distribution $PD(\theta)$. The sequence $X_k, k = 1, 2, \ldots$, corresponds to the size-biased permutation of $PD(\theta)$ and the representation through $U_k, k = 1, 2, \ldots$, is called the GEM representation after R. C. Griffiths, S. Engen and J. W. McCloskey.

Consider a population under the influence of mutation and selection. The role of mutation is to bring in new types of alleles and reduce the proportion of existing alleles, while the selection force favors certain genotypes and, thus, alters allele proportions. It is interesting to understand how the mutation and selection forces interact. The limiting procedure with $\theta$ approaching infinity is equivalent to letting the population size go to infinity. By the study of the behavior of $PD(\theta)$ for large $\theta$, one would hope to get a better picture of interaction between mutation and selection. For the overdominance model, where the heterozygote has advantage over homozygote, it is observed in [5] that, when both the mutation rate and the selection rate are scaled by large $\theta$, the model behaves the same as a neutral model. This was confirmed later by Joyce, Krone and Kurtz [9] through the study of the stationary distribution of the infinitely many alleles diffusion with heterozygote advantage. A critical growth rate $\theta^{3/2}$ is identified such that selection will not be detected if its rate grows more slowly than the critical rate.

Let $\xi_k, k = 1, \ldots$, be a sequence of i.i.d. random variables with common diffusive distribution $\nu$ on $[0, 1]$, that is, $\nu(\{x\}) = 0$ for every $x$ in $[0, 1]$. Set

$$\Pi = \sum_{k=1}^{\infty} P_k(\theta) \delta_{\xi_k}.$$

It is known that the law of $\Pi$ is the $Dirichlet(\nu)$ distribution, and is the stationary distribution of the Fleming–Viot process with mutation operator

$$Af(x) = \frac{\theta}{2}\int_0^1 (f(y) - f(x))\nu(dx).$$

Dawson and Feng [1, 2] studied the asymptotic behavior of $\Pi$ for large $\theta$ and established the large deviation principle (henceforth LDP) for the law of $\Pi$. It is worth noting that there are fundamental differences between $\Pi$ and $\mathbf{P}(\theta)$ even though their laws are both called Poisson–Dirichlet distribution in the literature. A detailed discussion is given in Section 4.

The main result of this article is the LDP for $PD(\theta)$ for large $\theta$. When $\theta$ approaches infinity, $P_k(\theta)$ converges to zero for every $k$. Since $\sum_{k=1}^\infty P_k = 1$, the allele proportions are evenly spread out for large $\theta$. We will see from the LDP that, at the exponential scale, the differences between different allele proportions are still significant.

Our first result is the LDP for the law of $P_1(\theta)$. This is then used to derive the LDP for finite marginal distributions of $\mathbf{P}(\theta)$, namely, the law of $(P_1(\theta), \ldots, P_n(\theta))$ for every $n$. These eventually lead to the establishment of the LDP for the law of $\mathbf{P}(\theta)$. All rate functions have explicit forms.

In Section 2 we review several general results on LDP, and formulate a comparison lemma. Some estimates on the Beta distribution are proved in Section 3. Our main LDP results for $PD(\theta)$ are formulated and proved in Section 4. In Section 5 the LDP result for $PD(\theta)$ is used to derive the LDP of the homozygosity and the $PD(\theta)$ with selection. Our result shows that a phase transition occurs with parameter given by the selection intensity. Let the selection be scaled by $\theta^\gamma$. Then for the selection to be detected at the large deviation scale, $\gamma$ has to be greater than or equal to 1. For the heterozygote advantage model, this provides a new proof of a conjecture in [5]. In the LDP setting the critical scale turns out to be $\theta$ instead of $\theta^{3/2}$ which was obtained in the case of Joyce, Krone and Kurtz [9].

The study of the behavior of $\mathbf{P}(\theta) = (P_1(\theta), P_2(\theta), \ldots)$ for large $\theta$ has a long history. In Waterson and Guess [17] $E[P_1(\theta)]$ was shown to be asymptotically $\log\theta/\theta$. Griffiths [6] obtained the explicit weak limit of $\theta\mathbf{P}(\theta)$ and a central limit theorem for the homozygosity. A more detailed description of these results and their relation to our results will be included in Section 4.

One may be able to generalize our result to the two-parameter Poisson–Dirichlet distribution studied in [12]. The residual allocation model now involves two parameters $\theta + \alpha > 0, 0 \leq \alpha < 1$, such that $U_k$ is a $Beta(1 - \alpha, \theta + k\alpha)$ random variable for each $k$. Since the mutation force becomes stronger with the introduction of $\alpha$, one expects the speed of convergence will be higher than that of $PD(\theta)$. For a more comprehensive discussion on $PD(\theta)$ and its two-parameter counterpart, we recommend [13, 14] and the references therein.

**2. Preliminaries.** We include several known results on LDP in this section. A comparison lemma will be formulated as a direct application of the Gärtner–Ellis theorem. All results will be stated in the form that is sufficient for our purposes. For the most general form, we refer to [3]. Let $E$ be a complete separable metric space with metric $\rho$.

DEFINITION 2.1.  A family of probability measures $\{Q_\varepsilon : \varepsilon > 0\}$ on $E$ is said to satisfy an LDP with speed $1/\varepsilon$ and rate function $I(\cdot)$ if, for any closed set $F$ and open set $G$ in $E$,

$$\limsup_{\varepsilon \to 0} \varepsilon \log Q_\varepsilon\{F\} \leq -\inf_{x \in F} I(x),$$

$$\liminf_{\varepsilon \to 0} \varepsilon \log Q_\varepsilon\{G\} \geq -\inf_{x \in G} I(x),$$

for any $c > 0$,   $\{x : I(x) \leq c\}$ is compact.

DEFINITION 2.2.  A family of probability measures $\{Q_\varepsilon : \varepsilon > 0\}$ is said to satisfy a partial LDP if, for every sequence $\varepsilon_n$ converging to zero, there is a subsequence $\varepsilon'_n$ such that the family $\{Q_{\varepsilon'_n} : \varepsilon'_n > 0\}$ satisfies an LDP with speed $1/\varepsilon'_n$ and rate function $I'$.

REMARK.  A partial LDP will become an LDP if all the rate functions $I'$ are the same. The following result is found in [15].

THEOREM 2.1 (Pukhalskii).  (i) *Assume that $\{Q_\varepsilon : \varepsilon > 0\}$ satisfies a partial LDP with speed $1/\varepsilon$ and for every $x$ in $E$,*

$$\lim_{\delta \to 0} \limsup_{\varepsilon \to 0} \varepsilon \log Q_\varepsilon\{\rho(y,x) \leq \delta\}$$
(1)
$$= \lim_{\delta \to 0} \liminf_{\varepsilon \to 0} \varepsilon \log Q_\varepsilon\{\rho(y,x) < \delta\} = -I(x).$$

*Then $\{Q_\varepsilon : \varepsilon > 0\}$ satisfies an LDP with speed $1/\varepsilon$ and rate function $I(\cdot)$.*
(ii) *If $E$ is compact, then the partial large deviation principle is automatically satisfied.*

THEOREM 2.2 (Varadhan).  *Assume that $\{Q_\varepsilon : \varepsilon > 0\}$ satisfy an LDP with speed $1/\varepsilon$ and a rate function $I(\cdot)$. Let $C_b(E)$ denote the set of bounded continuous functions on $E$. Then for any $\phi(x)$ in $C_b(E)$, one has*

$$\Lambda_\phi = \lim_{\varepsilon \to 0} \varepsilon \log E^{Q_\varepsilon}[e^{\phi(x)/\varepsilon}] = \sup_{x \in E}\{\phi(x) - I(x)\}. \tag{2}$$

THEOREM 2.3 (Gärtner–Ellis).  *Let $E = R = (-\infty, \infty)$. Assume that*

$$\Lambda(\lambda) = \lim_{\varepsilon \to 0} \varepsilon \log E^{Q_\varepsilon}[e^{\lambda x/\varepsilon}] < \infty \qquad \text{for all } \lambda,$$

and has the first-order derivative $\Lambda'(\lambda)$. Then $\{Q_\varepsilon : \varepsilon > 0\}$ satisfies an LDP with speed $1/\varepsilon$ and rate function

$$I(x) = \sup_{\lambda \in E}\{\lambda x - \Lambda(\lambda)\}.$$

As a direct application of Theorem 2.3, we get the following:

LEMMA 2.4 (Comparison lemma). *Let $E = (-\infty, \infty)$. Assume that $\{X_\varepsilon : \varepsilon > 0\}$, $\{Y_\varepsilon : \varepsilon > 0\}$, $\{Z_\varepsilon : \varepsilon > 0\}$ are three families of random variables on the same probability space with respective laws $\{Q_\varepsilon^1 : \varepsilon > 0\}$, $\{Q_\varepsilon^2 : \varepsilon > 0\}$, $\{Q_\varepsilon^3 : \varepsilon > 0\}$. If both $\{Q_\varepsilon^1 : \varepsilon > 0\}$ and $\{Q_\varepsilon^3 : \varepsilon > 0\}$ satisfy the assumptions in Theorem 2.3 with the same $\Lambda$, and with probability one*

$$X_\varepsilon \leq Y_\varepsilon \leq Z_\varepsilon,$$

*then $\{Q_\varepsilon^2 : \varepsilon > 0\}$ satisfies an LDP with speed $1/\varepsilon$ and rate function*

$$I(x) = \sup_{\lambda \in E}\{\lambda x - \Lambda(\lambda)\}.$$

REMARK. Both Theorem 2.3 and Lemma 2.4 hold if $E$ is only a closed subset of $R$.

**3. Some estimates for the beta distribution.** Let $U_1, U_2, \ldots$ be a sequence of i.i.d. random variables with common distribution $Beta(1, \theta)$. Let $E = [0, 1]$. Then we have the following:

LEMMA 3.1. *For any $n \geq 1$, let $Q_{n,\theta}$ be the law of $Z_n = \max\{U_1, \ldots, U_n\}$. Then the family $\{Q_{n,\theta} : \theta > 0\}$ satisfies an LDP on $E$ with speed $\theta$ and rate function*

$$(3) \qquad I(x) = \begin{cases} \log \dfrac{1}{1-x}, & x \in [0, 1), \\ \infty, & else. \end{cases}$$

PROOF. Let

$$\Lambda(\lambda) = \operatorname*{ess\,sup}_{y \in [0,1]}\{\lambda y + \log(1-y)\}$$

$$(4) \qquad = \begin{cases} \lambda - 1 - \log \lambda, & \lambda > 1, \\ 0, & else. \end{cases}$$

Then clearly $\Lambda(\lambda)$ is finite for all $\lambda$ and is differentiable. By direct calculation, we have

$$E[e^{\theta \lambda Z_n}] = \int_0^1 \exp\{\theta F_\theta(y)\}\,dy,$$

where

$$F_\theta(y) = \lambda y + \frac{\log n + \log \theta}{\theta}$$

(5)

$$+ \frac{n-1}{\theta}\log[1 - (1-y)^\theta] + \frac{\theta - 1}{\theta}\log(1-y).$$

Letting $\theta$ go to infinity, we get

$$\lim_{\theta \to \infty} \log\{E[e^{\theta \lambda Z_n}]\}^{1/\theta} = \Lambda(\lambda)$$

which, combined with Theorem 2.3 (with $\varepsilon = 1/\theta$), implies the lemma. □

LEMMA 3.2. *For any $k \geq 1$, let $n_k(\theta)$ denote the integer part of $\theta^k$. Then the family $\{Q_{n_k(\theta),\theta} : \theta > 0\}$ satisfies an LDP with speed $\theta$ and rate function $I(\cdot)$ defined in (3).*

PROOF. Choosing $n = n_k(\theta)$ in Lemma 3.1, we get

$$F_\theta(y) = \lambda y + \frac{(\log n_k(\theta) + \log \theta)}{\theta}$$

$$+ \frac{n_k(\theta) - 1}{\theta}\log[1 - (1-y)^\theta] + \frac{\theta - 1}{\theta}\log(1-y).$$

For any $\varepsilon$ in $(0, 1/2)$, and $\lambda \geq 0$, we have

$$\Lambda(\lambda) = \lim_{\theta \to \infty} \frac{1}{\theta} \log E[e^{\theta \lambda U_1}]$$

$$\leq \limsup_{\theta \to \infty} \frac{1}{\theta} \log E[e^{\theta \lambda Z_{n_k(\theta)}}]$$

$$\leq \max\left\{\lambda \varepsilon, \operatorname*{ess\,sup}_{y \geq \varepsilon}[\lambda y + \log(1-y)]\right\},$$

where the last inequality follows from the fact that, for $y$ in $[\varepsilon, 1]$, $\lim_{\theta \to \infty} \theta^l \times \log[1 - (1-y)^\theta] = 0$ for any $l \geq 1$. Letting $\varepsilon$ go to zero, it follows that

$$\Lambda(\lambda) = \lim_{\theta \to \infty} \frac{1}{\theta} \log E[e^{\theta \lambda Z_{n_k(\theta)}}].$$

For negative $\lambda$, we have

$$\limsup_{\theta \to \infty} \frac{1}{\theta} \log E[e^{\theta \lambda Z_{n_k(\theta)}}] \geq \lim_{\delta \to 0} \operatorname*{ess\,sup}_{y \geq \delta}\{\lambda y + \log(1-y)\} = 0 = \Lambda(\lambda).$$

The lemma follows from Theorem 2.3. □

LEMMA 3.3. *For any $n \geq 1$, let $W_n = (1-U_1)(1-U_2)\cdots(1-U_n)$. Then for any $\delta > 0$,*

(6) $$\limsup_{\theta \to \infty} \frac{1}{\theta} \log P\{W_{n_2(\theta)} \geq \delta\} = -\infty.$$

PROOF. By direct calculation,

$$P\{W_{n_2(\theta)} \geq \delta\} = P\left\{\theta \sum_{j=1}^{n_2(\theta)} \log(1-U_j) \geq \theta \log \delta\right\}$$

$$\leq e^{\theta \log 1/\delta}(E[e^{\theta \log(1-U_1)}])^{n_2(\theta)}$$

$$= e^{\theta \log 1/\delta}\left(\frac{1}{2}\right)^{n_2(\theta)}$$

$$= \exp\left[\theta \log \frac{1}{\delta} - (\theta^2 - 1)\log 2\right].$$

(6) follows by letting $\theta$ go to infinity. □

**4. LDP for the Poisson–Dirichlet distribution.** Let $\mathbf{X}(\theta) = (X_1, X_2, \ldots)$ be the GEM, that is,

(7) $\qquad X_1 = U_1, \qquad X_k = (1-U_1)\cdots(1-U_{k-1})U_k, \qquad k \geq 2,$

and

(8) $$\mathbf{P}(\theta) = (P_1(\theta), P_2(\theta), \ldots),$$

with $P_k(\theta)$ the $k$th largest component of $\mathbf{X}(\theta)$. The law of $\mathbf{P}(\theta)$ is thus $PD(\theta)$.

In this section we will establish the LDP for $PD(\theta)$ when $\theta$ becomes large. To help motivate this result, some earlier works on the asymptotic behavior of $PD(\theta)$ are included and their relations to our result are discussed.

4.1. *Scaling limits.* Recall that the parameter $\theta$ is the population mutation rate. In the infinite neutral allele models all mutations produce new alleles. It is thus reasonable to expect that the higher the mutation rate, the smaller the proportion of most frequent allele will be.

In [17], the exact expression and asymptotic expression were obtained for $E[P_1(\theta)]$. In particular, they showed that

(9) $$\lim_{\theta \to \infty} \frac{E[P_1(\theta)]}{(\log \theta/\theta)} = 1,$$

which implies that $\lim_{\theta \to \infty} P_k(\theta) = 0$ for each $k$. Since $\sum_{k=1}^{\infty} P_k(\theta) = 1$, it follows that the differences between the proportions $\{P_k(\theta) : k \geq 1\}$ become smaller when $\theta$ becomes large.

Griffiths [6] generalized the result in [17], and obtained expressions for the expectations and variances of $P_k(\theta)$ for each $k$. The moments of $P_k(\theta)$ for any $k \geq 1$ were given by the following:

$$
(10) \qquad E[P_k^n(\theta)] = \frac{\theta^k \Gamma(\theta)}{\Gamma(\theta + n)} \int_0^{\infty} u^{n-1} \frac{[J(u)]^{k-1}}{(k-1)!} e^{-u - \theta J(u)} \, du,
$$

where $J(u) = \int_u^{\infty} \frac{e^{-x}}{x} \, dx$. In particular, one has, for any $k \geq 1$,

$$
E[P_k(\theta)] = \int_0^{\infty} e^{-u} \left[ \frac{(\theta J(u))^{k-1}}{(k-1)!} e^{-\theta J(u)} \right] du \to 0 \qquad \text{as } \theta \to \infty.
$$

Perman [11] obtained generalizations of (10) to normalized jump sizes of subordinators. A scaling limit, to be described below, was also obtained in [6].

For each $r \geq 1$, let $\infty > Y_1 > Y_2 > \cdots > Y_r > -\infty$ have a joint distribution with density

$$
(11) \qquad \exp\{-(y_1 + \cdots + y_r) - e^{-y_r}\}.
$$

It is clear that the marginal density of $Y_k$ is

$$
(12) \qquad \frac{1}{(k-1)!} \exp\{-(ky + e^{-y})\}.
$$

For $k = 1, \ldots, r$, set

$$
(13) \qquad \beta(\theta) = \log \theta + \log \log \theta, \qquad Y_k(\theta) = \theta P_k(\theta) - \beta(\theta).
$$

THEOREM 4.1 (Griffiths). *For each $r \geq 1$, $(Y_1(\theta), \ldots, Y_r(\theta))$ converges weakly to $(Y_1, \ldots, Y_r)$ when $\theta$ goes to infinity.*

The result (9) can be viewed as a kind of law of large numbers and Theorem 4.1 as a "central limit" type theorem. This brings us naturally to the study of large deviations in the next subsection.

4.2. *Large deviations.* There are two different versions of the infinitely-many-neutral-alleles model: one is a special Fleming–Viot process with parent independent mutation operator with mutation rate $\theta$ and mutation probability $\nu$, and the other is an infinite-dimensional diffusion process with state space $\nabla$ and generator

$$
L = \frac{1}{2} \sum_{i,j=1}^{\infty} p_i(\delta_{ij} - p_j) \frac{\partial^2}{\partial p_i \, \partial p_j} - \frac{\theta}{2} \sum_{i=1}^{\infty} \frac{\partial}{\partial p_i}
$$

defined on an appropriate domain. The Fleming–Viot version is called the labeled version and the second version is called unlabeled. Fundamental differences exist between the two versions. For example, the labeled version does not have a transition density, while the unlabeled version does; the unlabeled version has one less eigenvalue than the labeled version (see [4]). But both models are reversible with respective reversible measures $Dirichlet(\nu)$ and $PD(\theta)$. Let $M_1([0,1])$ be the set of all probability measures on $[0,1]$. If we introduce the map $\Phi$ between $M_1([0,1])$ and the closure of $\nabla$ in $R^\infty$ such that $\Phi(\mu)$ is the descending sequence of masses of the atoms of $\mu$, then the unlabeled model is just the image of the labeled Fleming–Viot model under $\Phi$. Thus, many properties of the unlabeled version can be derived from the labeled one. Since the LDP for $Dirichlet(\nu)$ has been established in [1, 2], one would hope to get the LDP for $PD(\theta)$ from the LDP for $Dirichlet(\nu)$ through $\Phi$. Unfortunately, $\Phi$ is not continuous as the following example shows. Let $\mu_n = \frac{1}{n}\sum_{k=1}^n \delta_{k/n^2}$. Then $\mu_n$ converges weakly to $\delta_0$, while $\Phi(\mu_n) = (\frac{1}{n},\ldots,\frac{1}{n},0,0\ldots)$ converges to $(0,\ldots 0,\ldots)$ rather than $(1,0,\ldots) = \Phi(\delta_0)$. To use the contraction principle in large deviation theory, one has to prove some exponential approximation to $\Phi$ by a sequence of continuous maps. We choose to prove the LDP for $PD(\theta)$ directly.

Our first theorem gives the large deviations of $P_1(\theta)$.

LEMMA 4.2. *The family of the laws of $P_1(\theta)$ satisfies an LDP on $[0,1]$ with speed $\theta$ and rate function $I(\cdot)$ [given by (3)].*

PROOF. Let $\hat{P}_1(\theta) = \max\{X_1,\ldots,X_{n_2(\theta)}\}$. Then clearly $P_1(\theta) \geq \hat{P}_1(\theta)$. By Lemma 3.3, for any $\delta > 0$, one has

$$\limsup_{\theta \to \infty} \frac{1}{\theta} \log P\{P_1(\theta) - \hat{P}_1(\theta) > \delta\} \leq \limsup_{\theta \to \infty} \frac{1}{\theta} \log P\{W_{n_2(\theta)} > \delta\} = -\infty.$$

In other words, $P_1(\theta)$ and $\hat{P}_1(\theta)$ are exponentially equivalent and, thus, have the same LDPs. By definition, we have

$$U_1 = X_1 \leq \hat{P}_1(\theta) \leq Z_{n_2(\theta)}.$$

Applying Lemmas 3.1, 3.2 and 2.4, we conclude that the law of $\hat{P}_1(\theta)$ satisfies an LDP with speed $\theta$ and rate function $I(\cdot)$. □

Let

(14) $$\nabla_n = \left\{(p_1,\ldots,p_n) : 0 \leq p_n \leq \cdots \leq p_1, \sum_{k=1}^n p_k \leq 1\right\}$$

and

(15) $\Xi_{n,\theta}$ is the law of $(P_1(\theta),\ldots,P_n(\theta))$ on space $\nabla_n$, $\quad n \geq 1$.

THEOREM 4.3. *For fixed $n \geq 2$, the family $\{\Xi_{n,\theta} : \theta > 0\}$ satisfies an LDP with speed $\theta$ and rate function*

$$(16) \quad S_n(p_1,\ldots,p_n) = \begin{cases} \log \dfrac{1}{1 - \sum_{k=1}^n p_k}, & (p_1,\ldots,p_n) \in \nabla_n, \sum_{k=1}^n p_k < 1, \\ \infty, & \text{else.} \end{cases}$$

PROOF. Since $\nabla_n$ is compact, by Theorem 2.1(ii), the family $\{\Xi_{n,\theta} : \theta > 0\}$ satisfies a partial LDP. Let $g_1^\theta$ denote the density function of $P_1(\theta)$. Then for any $p \in (0,1)$,

$$(17) \quad g_1^\theta(p) p (1-p)^{1-\theta} = \theta \int_0^{(p/(1-p)) \wedge 1} g_1^\theta(x)\, dx,$$

and the joint density function $g_n^\theta$ of $(P_1(\theta),\ldots,P_n(\theta))$ obtained in Watterson [16] is given by

$$(18) \quad g_n^\theta(p_1,\ldots,p_n) = \frac{\theta^{n-1}(1 - \sum_{k=1}^{n-1} p_k)^{\theta-2}}{p_1 \cdots p_{n-1}} g_1^\theta\left(\frac{p_n}{1 - \sum_{k=1}^{n-1} p_k}\right),$$

for $(p_1,\ldots,p_n) \in \nabla_n^\circ = \{(p_1,\ldots,p_n) \in \nabla_n : 0 < p_n < \cdots < p_1 < 1, \sum_{k=1}^n p_k < 1\}$, and is zero otherwise. In other words, for any fixed $(p_1,\ldots,p_n) \in \nabla_n^\circ$, we have

$$(19) \quad g_n^\theta(p_1,\ldots,p_n) = \frac{\theta^n (1 - \sum_{k=1}^n p_k)^{\theta-1}}{p_1 \cdots p_n} \int_0^{(p_n/(1-\sum_{k=1}^n p_k)) \wedge 1} g_1^\theta(u)\, du.$$

Clearly, $\nabla_n$ is the closure of $\nabla_n^\circ$. Now for any $(p_1,\ldots,p_n) \in \nabla_n$, let

$$V((p_1,\ldots,p_n); \delta) = \{(q_1,\ldots,q_n) \in \nabla_n : |q_k - p_k| < \delta, k = 1,\ldots,n\},$$
$$U((p_1,\ldots,p_n); \delta) = \{(q_1,\ldots,q_n) \in \nabla_n : |q_k - p_k| \leq \delta, k = 1,\ldots,n\}.$$

Then the family $\{V((p_1,\ldots,p_n); \delta) : \delta > 0, (p_1,\ldots,p_n) \in \nabla_n\}$ is a base for the topology of $\nabla_n$. Now assume that $p_n > 0$ and $\delta$ is smaller that $p_n$. By (19), we have that, for any $(q_1,\ldots,q_n)$ in $V((p_1,\ldots,p_n),\delta)$,

$$g_n^\theta(q_1,\ldots,q_n) \leq \frac{\theta^n (1 - \sum_{k=1}^n (p_k - \delta))^{\theta-1}}{(p_1 - \delta) \cdots (p_n - \delta)},$$

which implies

$$(20) \quad \limsup_{\theta \to \infty} \frac{1}{\theta} \log \Xi_{n,\theta}\{U((p_1,\ldots,p_n); \delta)\} \leq -\log \frac{1}{1 - \sum_{k=1}^n (p_k - \delta)}.$$

Letting $\delta$ go to zero, we get

$$(21) \quad \limsup_{\delta \to 0} \limsup_{\theta} \frac{1}{\theta} \log \Xi_{n,\theta}\{U((p_1,\ldots,p_n); \delta)\} \leq -S_n(p_1,\ldots,p_n).$$

Next we turn to lower bound. First noting that, if $\sum_{k=1}^{n} p_k = 1$, the lower bound is trivially true since $S_n(p_1, \ldots, p_n) = \infty$. Hence, we assume that $\sum_{k=1}^{n} p_k < 1$, $p_n > 0$. We also assume that $0 < \delta < (1 - \sum_{k=1}^{n} p_k)/n$. Using (19) again, one has that, for any $(q_1, \ldots, q_n)$ in $V((p_1, \ldots, p_n), \delta) \cap \nabla_n^\circ$,

$$g_n^\theta(q_1, \ldots, q_n) \geq \theta^n \left(1 - \sum_{k=1}^{n}(p_k + \delta)\right)^{\theta-1} \int_0^{((p_n-\delta)/(1-\sum_{k=1}^{n}(p_k-\delta)))\wedge 1} g_1^\theta(u) \, du,$$

which implies

$$\liminf_\theta \frac{1}{\theta} \log \Xi_{n,\theta}\{V((p_1,\ldots,p_n);\delta)\}$$
$$\geq -\log \frac{1}{1 - \sum_{k=1}^{n}(p_k+\delta)} - \inf\left\{I(p) : p < \left(\frac{(p_n - \delta)}{(1 - \sum_{k=1}^{n}(p_k - \delta))}\right) \wedge 1\right\}$$
$$= -\log \frac{1}{1 - \sum_{k=1}^{n}(p_k+\delta)},$$

where in the second line we used the LDP of the law $P_1(\theta)$ obtained in Lemma 4.2. Letting $\delta$ go to zero, we get

$$(22) \quad \liminf_{\delta \to 0} \liminf_{\theta \to \infty} \frac{1}{\theta} \log \Xi_{n,\theta}\{V((p_1,\ldots,p_n);\delta)\} \geq -S_n(p_1,\ldots,p_n).$$

Finally, we turn to the case when there is $1 \leq k \leq n$ such that $p_i > 0$ for $i = 1, \ldots k$ and $p_i = 0$ for $i \geq k+1$. Because of lower semi-continuity of all rate functions in the partial LDP and the continuity of $S_n(p_1, \ldots, p_n)$, (22) holds in this case. On the other hand, noting that $S_n(p_1, \ldots, p_n) = S_k(p_1, \ldots, p_k)$ and $\Xi_{n,\theta}\{U((p_1, \ldots, p_n); \delta)\} \leq \Xi_{k,\theta}\{U((p_1, \ldots, p_k); \delta)\}$, it follows that the upper bound also holds. By Theorem 2.1(i), (21) and (22) combined with the partial LDP imply the result. □

COROLLARY 4.1. *For $k \geq 2$, the family of the laws of $P_k(\theta)$ satisfies an LDP on $[0,1]$ with speed $\theta$ and rate function*

$$(23) \quad I_k(x) = \begin{cases} \log \dfrac{1}{1-kx}, & x \in [0, 1/k], \\ \infty, & else. \end{cases}$$

PROOF. For any $k \geq 2$, define

$$\phi_k : \nabla_k \longrightarrow [0,1], \qquad (p_1, p_2, \ldots, p_k) \to p_k.$$

Clearly, $\phi_k$ is continuous, and Theorem 4.3 combined with the contraction principle implies that the law of $P_k(\theta)$ satisfies an LDP on $[0,1]$ with speed $\theta$ and rate function

$$I'(p) = \inf\{S_k(p_1, \ldots, p_k) : p_1 \geq \cdots \geq p_k = p\}.$$

For $p > 1/k$, the infimum is over empty set and is thus infinity. For $p$ in $[0, 1/k]$, the infimum is achieved at the point $p_1 = p_2 = \cdots = p_k = p$. Hence, $I'(p) = I_k(p)$ and the result follows. $\square$

The result of this corollary indicates that, for any $k \geq 1$, the law of $kP_k(\theta)$ has the same LDP as the law of $P_1(\theta)$. More precisely, for any $p \in [0,1]$, one has

$$\lim_{\delta \to 0} \lim_{\theta \to \infty} \frac{P\{|P_1(\theta) - p| \leq \delta\}}{P\{|P_k(\theta) - p/k| \leq \delta\}} = 1. \tag{24}$$

Hence, when $\theta$ becomes large, $P_k(\theta)$ behaves like $\frac{1}{k}P_1(\theta)$ at the large deviation scale. In other words, under a large deviation, the proportion of the most likely alleles is $k$ times of the proportion of the $k$th most likely alleles.

This relation is also reflected somewhat in Theorem 4.1 and (12). We illustrate this through the following nonrigorous derivation with $\beta(\theta)$ defined in (13):

$$\begin{aligned} P\{P_1(\theta) \in dx\} &= P\{Y_1(\theta) + \beta(\theta) \in \theta\, dx\} \\ &\approx P\{Y_1 + \beta(\theta) \in \theta\, dx\} \\ &\approx \exp\left\{-\theta\left[x + \frac{\beta(\theta)}{\theta} + \frac{1}{\theta}e^{-\theta(x+\beta(\theta)/\theta)}\right]\right\} dx \end{aligned} \tag{25}$$

and

$$\begin{aligned} P\{P_k(\theta) \in dy\} &= P\{Y_k(\theta) + \beta(\theta) \in \theta\, dy\} \\ &\approx P\{k(Y_k + \beta(\theta)) \in \theta\, d(ky)\} \\ &\approx \exp\left\{-\theta\left[x + \frac{\beta(\theta)}{k\theta} + \frac{1}{k\theta}e^{-\theta/k(x+\beta(\theta)/\theta)}\right]\right\} dx, \qquad x = ky. \end{aligned} \tag{26}$$

Comparing the last terms in (25) and (26), we can see that at the exponential scale $kP_k(\theta)$ is like $P_1(\theta)$.

Now we turn to the LDP of $PD(\theta)$. Let

$$\bar{\nabla} = \left\{(p_1, p_2, \ldots) : p_1 \geq p_2 \geq \cdots \geq 0, \sum_{k=1}^{\infty} p_k \leq 1\right\} \tag{27}$$

be the closure of $\nabla$ equipped with the subspace topology of $R^\infty$. Let

$$\Xi_\theta \text{ be the law of } \mathbf{P}(\theta) \text{ on space } \bar{\nabla}. \tag{28}$$

THEOREM 4.4.  *The family $\{\Xi_\theta : \theta > 0\}$ satisfies an LDP with speed $\theta$ and rate function*

$$
(29) \quad S(\mathbf{p}) = \begin{cases} \log \dfrac{1}{1 - \sum_{k=1}^{\infty} p_k}, & (p_1, p_2, \ldots) \in \bar{\nabla}, \sum_{k=1}^{\infty} p_k < 1, \\ \infty, & else. \end{cases}
$$

PROOF. Because $\bar{\nabla}$ is compact, by Theorem 2.1, it suffices to verify (1) for the family $\{\Xi_\theta : \theta > 0\}$. The topology on $\bar{\nabla}$ can be generated by the following metric:

$$d(\mathbf{p}, \mathbf{q}) = \sum_{k=1}^{\infty} \frac{|p_k - q_k|}{2^k},$$

where $\mathbf{p} = (p_1, p_2, \ldots), \mathbf{q} = (q_1, q_2, \ldots)$. For any fixed $\delta > 0$, let $B(\mathbf{p}, \delta)$ and $\bar{B}(\mathbf{p}, \delta)$ denote the respective open and closed balls centered at $\mathbf{p}$ with radius $\delta > 0$. Set $n_\delta = 1 + [\log_2(1/\delta)]$, where $[x]$ denotes the integer part of $x$. Set

$$V_{n_\delta}(\mathbf{p}; \delta/2) = \{(q_1, q_2, \ldots) \in \bar{\nabla} : |q_k - p_k| < \delta/2, k = 1, \ldots, n_\delta\}.$$

Then we have

$$V_{n_\delta}(\mathbf{p}; \delta/2) \subset B(\mathbf{p}, \delta)$$

By Theorem 4.3 and the fact that

$$\Xi_\theta\{V_{n_\delta}(\mathbf{p}; \delta/2)\} = \Xi_{n_\delta, \theta}\{V((p_1, \ldots, p_{n_\delta}); \delta/2)\},$$

we get that

$$
\begin{aligned}
& \liminf_{\theta \to \infty} \frac{1}{\theta} \log \Xi_\theta\{B(\mathbf{p}, \delta)\} \\
(30) \quad & \geq \liminf_{\theta \to \infty} \frac{1}{\theta} \log \Xi_{n_\delta, \theta}\{V((p_1, \ldots, p_{n_\delta}); \delta/2)\} \\
& \geq -S_{n_\delta}(p_1, \ldots, p_{n_\delta}) \geq -S(\mathbf{p}).
\end{aligned}
$$

On the other hand, for any fixed $n \geq 1, \delta_1 > 0$, let

$$U_n(\mathbf{p}; \delta_1) = \{(q_1, q_2, \ldots) \in \bar{\nabla} : |q_k - p_k| \leq \delta_1, k = 1, \ldots, n\}.$$

Then we have

$$\Xi_\theta\{U_n(\mathbf{p}; \delta_1)\} = \Xi_{n, \theta}\{U((p_1, \ldots, p_n); \delta_1)\},$$

and, for $\delta$ small enough,

$$\bar{B}(\mathbf{p}, \delta) \subset U_n(\mathbf{p}; \delta_1),$$

which implies that

$$\lim_{\delta \to 0} \limsup_{\theta \to \infty} \frac{1}{\theta} \log \Xi_\theta \{\bar{B}(\mathbf{p}, \delta)\}$$

(31)
$$\leq \limsup_{\theta \to \infty} \frac{1}{\theta} \log \Xi_{n,\theta}\{U((p_1, \ldots, p_n), \delta_1)\}$$

$$\leq -\inf\{S_n(q_1, \ldots, q_n) : (q_1, \ldots, q_n) \in U((p_1, \ldots, p_n), \delta_1)\}.$$

Letting $\delta_1$ go to zero, and then $n$ go to infinity, we get

(32) $$\lim_{\delta \to 0} \limsup_{\theta \to \infty} \frac{1}{\theta} \log \Xi_\theta\{\bar{B}(\mathbf{p}, \delta)\} \leq -S(\mathbf{p}),$$

which combined with (30) implies the result. □

REMARK. Note that the effective domain is

$$\{\mathbf{p} \in \bar{\nabla} : S(\mathbf{p}) < \infty\} = \left\{\mathbf{p} \in \bar{\nabla} : \sum_{k=1}^{\infty} p_k < 1\right\}$$

and

$$\lim_{\delta \to 0} \inf_{\{\mathbf{p} : |\sum_{k=1}^{\infty} p_k - 1| \leq \delta\}} S(\mathbf{p}) = \infty.$$

On the other hand, since

$$\Xi_\theta\left\{\mathbf{p} \in \bar{\nabla} : \left|\sum_{k=1}^{\infty} p_k - 1\right| < \delta\right\} = \Xi_\theta\left\{\mathbf{p} \in \bar{\nabla} : \left|\sum_{k=1}^{\infty} p_k - 1\right| \leq \delta\right\} = 1,$$

one has

$$\lim_{\delta \to 0} \liminf_{\theta \to \infty} \frac{1}{\theta} \log \Xi_\theta\left\{\mathbf{p} \in \bar{\nabla} : \left|\sum_{k=1}^{\infty} p_k - 1\right| < \delta\right\}$$

$$= \lim_{\delta \to 0} \limsup_{\theta \to \infty} \frac{1}{\theta} \log \Xi_\theta\left\{\mathbf{p} \in \bar{\nabla} : \left|\sum_{k=1}^{\infty} p_k - 1\right| \leq \delta\right\} = 0.$$

This might at first sight appear to be a contradiction. However, since the function $\sum_{k=1}^{\infty} p_k$ is not continuous on $\bar{\nabla}$, the set $\{\mathbf{p} : |\sum_{k=1}^{\infty} p_k - 1| \leq \delta\}$ is not closed and there is no inconsistency.

**5. Applications.** In this section we will discuss two applications of Theorem 4.4. The first one is the LDP for the homozygosity.

A random sample of size $m > 1$ is selected from a population whose allelic types have distribution $PD(\theta)$. The probability that all samples are of the same type is called the $m$th order population homozygosity and is given by

(33) $$H_m(\theta) = \sum_{i=1}^{\infty} P_i^m(\theta) = \sum_{i=1}^{\infty} X_i^m.$$

Since $H_m(\theta) \leq P_1^{m-1}(\theta)$, it follows that $H_m(\theta)$ converges to zero as $\theta$ approaches to infinity. In [8] it is shown that $\frac{\theta^{m-1}}{\Gamma(m)} H_m(\theta)$ converges to 1 in probability, that is, $H_m(\theta)$ goes to zero at a magnitude of $\frac{\Gamma(m)}{\theta^{m-1}}$. Our next theorem describes the large deviations of $H_m(\theta)$ from zero.

THEOREM 5.1. *The law of $H_m(\theta)$ for $m > 1$ satisfies an LDP on $[0,1]$ with speed $\theta$ and rate function $I(y^{1/m})$, where $I(\cdot)$ is given by (3).*

PROOF. For $m > 1$ the map

$$\phi_m : \bar{\nabla} \longrightarrow [0,1], \qquad \mathbf{p} \to \sum_{k=1}^{\infty} p_k^m$$

is continuous. By Theorem 4.4 and the contraction principle, it follows that the law of $H_m(\theta)$ satisfies an LDP with speed $\theta$ and rate function

$$\bar{I}(y) = \inf\{S(\mathbf{p}) : \mathbf{p} \in \bar{\nabla}, \phi_m(\mathbf{p}) = y\}.$$

Since for any $\mathbf{p}$ in $\bar{\nabla}$, we have

$$\sum_{k=1}^{\infty} p_k \geq (\phi_m(\mathbf{p}))^{1/m} = y^{1/m},$$

it follows that $S(\mathbf{p}) \geq I(y^{1/m})$ and, thus, $\bar{I}(y) \geq I(y^{1/m})$. On the other hand, choosing $\mathbf{p} = (y^{1/m}, 0, \ldots)$, one gets that $\bar{I}(y) \leq I(y^{1/m})$. Hence, $\bar{I}(y) = I(y^{1/m})$, and the result follows. $\square$

The study of fluctuations of homozygosity goes back to Griffiths [6]. It was shown in [6] that

$$(34) \qquad \frac{\theta^{3/2}}{\sqrt{2}} [H_2(\theta) - E(H_2(\theta))] \to \mathcal{Z},$$

where $\mathcal{Z}$ is the standard normal random variable.

REMARK. It is interesting to note that the relation between the large deviation Theorem 5.1 and the "central limit theorem" (34) is qualitatively different from the corresponding relation in the classical case of partial sums of i.i.d. random variables. In the latter case the speed in the large deviation result is the same as that for the normal approximation and only the rate functions are different. In contrast, for the case of $H_2(\theta)$, the speed in the large deviation result is $\theta$ whereas for the normal approximation, $E(H_2(\theta)) + \frac{\sqrt{2}\mathcal{Z}}{\theta^{3/2}}$ it would be $\theta^3$.

Joyce, Krone and Kurtz [8] obtained the following generalization of (34) to the $m$th order homozygosity:

$$\sqrt{\theta}\left(\frac{\theta^{m-1}}{\Gamma(m)}H_m(\theta) - 1\right) \to \mathcal{Z}(m), \tag{35}$$

where $\mathcal{Z}(m)$ is a normal random variable with mean zero and variance $\frac{\Gamma(2m)}{\Gamma^2(m)} - m^2$. One can rewrite (35) as

$$\frac{\theta^{m-1/2}}{\Gamma(m)}(H_m(\theta) - E[H_m(\theta)]) \to \mathcal{Z}(m), \tag{36}$$

which includes (34) as a special case. Thus, we have two different laws of large numbers and two different central limit theorems: the convergence of $H_m(\theta)$ to zero and the fluctuations around the mean, and the convergence of $\frac{\theta^{m-1}}{\Gamma(m)}H_m(\theta)$ to one and the associated fluctuations. One can easily go from one to the other by a simple algebraic transformation.

The LDP obtained in Theorem 5.1 is associated with the convergence of $H_m(\theta)$ to zero. It is thus natural to expect obtaining an LDP result for convergence of $\frac{\theta^{m-1}}{\Gamma(m)}H_m(\theta)$ to one from Theorem 5.1. Unfortunately, we are unable to verify this, but have the following partial information about the possible candidate for the LDP speed if there is one.

Assume that an LDP with speed $\alpha(\theta)$ and rate function $I(\cdot)$ holds for the convergence of $\frac{\theta^{m-1}}{\Gamma(m)}H_m(\theta)$ to one. Then for any constant $c > 0$,

$$\begin{aligned}
P\left\{\frac{\theta^{m-1}}{\Gamma(m)}H_m(\theta) \geq 1 + c\right\} &\geq P\left\{\frac{\theta^{m-1}}{\Gamma(m)}X_1^m \geq 1 + c\right\} \\
&= P\left\{X_1 \geq \left(\frac{\Gamma(m)(1+c)}{\theta^{m-1}}\right)^{1/m}\right\} \\
&= \left[\left(1 - \frac{(\Gamma(m)(1+c))^{1/m}}{\theta^{(m-1)/m}}\right)^{\theta^{(m-1)/m}}\right]^{\theta^{1/m}},
\end{aligned} \tag{37}$$

which implies that

$$\inf_{x \geq 1+c} I(x) = 0 \qquad \text{if } \lim_{\theta \to \infty}\frac{\alpha(\theta)}{\theta^{1/m}} = \infty. \tag{38}$$

Since $c$ is arbitrary, $I(\cdot)$ is zero over a sequence that goes to infinity, which contradicts the fact that $\{x : I(x) \leq M\}$ is compact for every positive $M$. Hence, the LDP speed cannot grow faster than $\theta^{1/m}$.

Our second application involves the Poisson–Dirichlet distribution with selection.

Let $C_b(\bar{\nabla})$ be the set of all bounded continuous functions on $\bar{\nabla}$. Assume that $\alpha(\theta)$ satisfy either $\lim_{\theta \to \infty} \frac{\alpha(\theta)}{\theta} = 0$ or $\lim_{\theta \to \infty} \frac{\alpha(\theta)}{\theta} = c$ for some $c > 0$. For every $H$ in $C_b(\bar{\nabla})$, define a new probability $\Xi_{\alpha,\theta}^H$ on $\bar{\nabla}$ as

$$(39) \qquad \Xi_{\alpha,\theta}^H(d\mathbf{p}) = \frac{e^{\alpha(\theta)H(\mathbf{p})}}{E^{\Xi_\theta}[e^{\alpha(\theta)H(\mathbf{q})}]} \Xi_\theta(d\mathbf{p}).$$

Then we have the following:

THEOREM 5.2. *The family $\{\Xi_{\alpha,\theta}^H\}$ satisfies an LDP with speed $\theta$ and rate function*

$$(40) \qquad S_\alpha(\mathbf{p}) = \begin{cases} S(\mathbf{p}), & \text{if } \lim_{\theta \to \infty} \frac{\alpha(\theta)}{\theta} = 0, \\ S^c(\mathbf{p}), & \text{if } \lim_{\theta \to \infty} \frac{\alpha(\theta)}{\theta} = c > 0, \end{cases}$$

*where*

$$(41) \qquad S^c(\mathbf{p}) = \sup_{\mathbf{q}}\{cH(\mathbf{q}) - S(\mathbf{q})\} - (cH(\mathbf{p}) - S(\mathbf{p})).$$

PROOF. By Theorem 2.2 and Theorem 4.4,

$$\lim_{\theta \to \infty} \frac{1}{\theta} \log E^{\Xi_\theta}[e^{\alpha(\theta)H(\mathbf{p})}]$$

$$= \lim_{\theta \to \infty} \frac{1}{\theta} \log E^{\Xi_\theta}[e^{\theta(\alpha(\theta)/\theta)H(\mathbf{p})}]$$

$$= \begin{cases} 0, & \text{if } \lim_{\theta \to \infty} \frac{\alpha(\theta)}{\theta} = 0, \\ \sup_{\mathbf{q}}\{cH(\mathbf{q}) - S(\mathbf{q})\}, & \text{if } \lim_{\theta \to \infty} \frac{\alpha(\theta)}{\theta} = c > 0. \end{cases}$$

This, combined with the continuity of $H$, implies that, for any $\mathbf{p}$ in $\bar{\nabla}$,

$$\liminf_{\delta \to 0} \liminf_{\theta \to \infty} \frac{1}{\theta} \log \Xi_{\alpha,\theta}^H\{d(\mathbf{p},\mathbf{q}) < \delta\}$$

$$\geq \liminf_{\delta \to 0} \liminf_{\theta \to \infty} \left\{ \frac{\alpha(\theta)}{\theta}(H(\mathbf{p}) - \delta') + \frac{1}{\theta} \log \Xi_\theta\{d(\mathbf{p},\mathbf{q}) < \delta\} \right\}$$

$$- \begin{cases} 0, & \text{if } \lim_{\theta \to \infty} \frac{\alpha(\theta)}{\theta} = 0, \\ \sup_{\mathbf{q}}\{cH(\mathbf{q}) - S(\mathbf{q})\}, & \text{if } \lim_{\theta \to \infty} \frac{\alpha(\theta)}{\theta} = c > 0, \end{cases}$$

$$\geq \begin{cases} -S(\mathbf{p}), & \text{if } \lim_{\theta \to \infty} \frac{\alpha(\theta)}{\theta} = 0, \\ -S^c(\mathbf{p}), & \text{if } \alpha(\theta) = c\theta > 0, \end{cases}$$

where $\delta'$ converges to zero as $\delta$ goes to zero. Similarly, we have

$$\limsup_{\delta \to 0} \limsup_{\theta \to \infty} \frac{1}{\theta} \log \Xi_{\alpha,\theta}^H \{d(\mathbf{p}, \mathbf{q}) \leq \delta\}$$

$$\leq \limsup_{\delta \to 0} \limsup_{\theta \to \infty} \left\{ \frac{\alpha(\theta)}{\theta}(H(\mathbf{p}) + \delta') + \frac{1}{\theta} \log \Xi_\theta \{d(\mathbf{p}, \mathbf{q}) \leq \delta\} \right\}$$

$$- \begin{cases} 0, & \text{if } \lim_{\theta \to \infty} \frac{\alpha(\theta)}{\theta} = 0, \\ \sup_{\mathbf{q}} \{cH(\mathbf{q}) - S(\mathbf{q})\}, & \text{if } \lim_{\theta \to \infty} \frac{\alpha(\theta)}{\theta} = c > 0, \end{cases}$$

$$\leq \begin{cases} -S(\mathbf{p}), & \text{if } \lim_{\theta \to \infty} \frac{\alpha(\theta)}{\theta} = 0, \\ -S^c(\mathbf{p}), & \text{if } \alpha(\theta) = c\theta > 0. \end{cases}$$

Since $\bar{\nabla}$ is compact, the result follows by an application of Theorem 2.1. $\square$

THEOREM 5.3. *Assume that*

$$\lim_{\theta \to \infty} \frac{\alpha(\theta)}{\theta} = \infty \tag{42}$$

*and the maximum of $H$ is achieved at a single point $\mathbf{p}_0$. Then the family $\{\Xi_{\alpha,\theta}^H\}$ satisfies an LDP with speed $\theta$ and rate function*

$$S^\infty(\mathbf{p}) = \begin{cases} 0, & \text{if } \mathbf{p} = \mathbf{p}_0, \\ \infty, & \text{else}. \end{cases} \tag{43}$$

PROOF. Without loss of generality, we assume that $\sup_{\mathbf{p} \in \bar{\nabla}} H(\mathbf{p}) = 0$. Otherwise we can multiply both the numerator and the denominator by $e^{-\alpha(\theta)H(\mathbf{p}_0)}$ in the definition of $\Xi_{\alpha,\theta}^H$.

For any $\mathbf{p} \neq \mathbf{p}_0$, choose $\delta$ small enough such that

$$d_1 = \sup_{d(\mathbf{p},\mathbf{q}) \leq \delta} H(\mathbf{q}) < 2d_2 = 2 \inf_{d(\mathbf{p}_0,\mathbf{q}) \leq \delta} H(\mathbf{q}) < 0.$$

Then by direct calculation, we get

$$\limsup_{\theta \to \infty} \frac{1}{\theta} \log \Xi_{\alpha,\theta}^H \{d(\mathbf{p}, \mathbf{q}) \leq \delta\}$$

$$= \limsup_{\theta \to \infty} \frac{1}{\theta} \log \frac{\int_{\{d(\mathbf{p},\mathbf{q}) \leq \delta\}} e^{\alpha(\theta)H(\mathbf{q})} \Xi_\theta(d\mathbf{q})}{E^{\Xi_\theta}[e^{\alpha(\theta)H(\mathbf{q})}]}$$

$$\leq \limsup_{\theta \to \infty} \frac{1}{\theta} \log \left[ e^{\alpha(\theta)(d_1 - d_2)} \frac{\Xi_\theta \{d(\mathbf{p}, \mathbf{q}) \leq \delta\}}{\Xi_\theta \{d(\mathbf{p}_0, \mathbf{q}) \leq \delta\}} \right]$$

$$= -\infty$$

and

$$\liminf_{\theta\to\infty}\frac{1}{\theta}\log\Xi_{\alpha,\theta}^H\{d(\mathbf{p}_0,\mathbf{q})<\delta\}$$

$$=\liminf_{\theta\to\infty}\frac{1}{\theta}\log\frac{\int_{\{d(\mathbf{p}_0,\mathbf{q})<\delta\}}e^{\alpha(\theta)H(\mathbf{q})}\Xi_\theta(d\mathbf{q})}{E^{\Xi_\theta}[e^{\alpha(\theta)H(\mathbf{q})}]}$$

$$=\liminf_{\theta\to\infty}\frac{1}{\theta}\log\left[1-\frac{\int_{\{d(\mathbf{p}_0,\mathbf{q})\geq\delta\}}e^{\alpha(\theta)H(\mathbf{q})}\Xi_\theta(d\mathbf{q})}{E^{\Xi_\theta}[e^{\alpha(\theta)H(\mathbf{q})}]}\right]$$

$$\geq\liminf_{\theta\to\infty}\frac{1}{\theta}\log\left[1-\frac{\int_{\{d(\mathbf{p}_0,\mathbf{q})\geq\delta\}}e^{\alpha(\theta)H(\mathbf{q})}\Xi_\theta(d\mathbf{q})}{\int_{\{d(\mathbf{p}_0,\mathbf{q})\leq\delta_1\}}e^{\alpha(\theta)H(\mathbf{q})}\Xi_\theta(d\mathbf{q})}\right]$$

$$\geq\liminf_{\theta\to\infty}\frac{1}{\theta}\log\left[1-\frac{\exp\{\alpha(\theta)[a(\delta)-b(\delta_1)]\}}{\Xi_\theta\{d(\mathbf{p}_0,\mathbf{q})\leq\delta\}}\right]$$

$$=0 \quad\text{by choosing small enough }\delta_1,$$

where

$$a(\delta)=\sup_{d(\mathbf{p}_0,\mathbf{q})\geq\delta}H(\mathbf{q}),\; b(\delta_1)=\inf_{d(\mathbf{p}_0,\mathbf{q})\leq\delta_1}H(\mathbf{q}).$$

This, combined with the compactness of $\bar{\nabla}$ and Theorem 2.1, implies the result. □

In [5], simulations were done for several models to study the role of population size in population genetical models of molecular evolution. One of the models is an infinite-alleles model with selective overdominance or heterozygote advantage. It was observed and conjectured that, if the selection intensity and the mutation rate get large at the same speed, the behavior looks like that of a neutral model. A rigorous proof of this conjecture was included in [9]. Using our notation with $\phi_m(\mathbf{p})=\sum_{k=1}^\infty p_k^m$, their result can be stated as follows.

THEOREM 5.4 (Joyce, Krone and Kurtz). *Choosing $\alpha(\theta)=c\theta^{3/2+\gamma}$ and $H(\mathbf{p})=-\phi_2(\mathbf{p})$ in (39), then, under $\Xi_\theta$, as $\theta\to\infty$,*

$$(44)\qquad \frac{e^{\alpha(\theta)H(\mathbf{p})}}{E^{\Xi_\theta}[e^{\alpha(\theta)H(\mathbf{p})}]}\Rightarrow\begin{cases}1,&\text{if }\gamma<0,\\ \exp(c\mathcal{Z}(2)-c^2),&\text{if }\gamma=0,\\ 0,&\text{if }\gamma>0,\end{cases}$$

*where $\Rightarrow$ denotes the weak convergence and $\mathcal{Z}(2)$ is a normal random variable with mean zero and variance 2.*

Now choosing $H(\mathbf{p})=-\phi_2(\mathbf{p})$ in Theorem 5.2 and Theorem 5.3, the next corollary gives an alternate proof of Gilliespie's conjecture.

COROLLARY 5.1. *The family $\Xi_{\alpha,\theta}^H$ satisfies an LDP with speed $\theta$ and rate function*

$$
(45) \quad S_\alpha(\mathbf{p}) = \begin{cases} S(\mathbf{p}), & \text{if } \lim_{\theta \to \infty} \frac{\alpha(\theta)}{\theta} = 0, \\ S^c(\mathbf{p}), & \text{if } \lim_{\theta \to \infty} \frac{\alpha(\theta)}{\theta} = c > 0, \\ S^\infty(\mathbf{p}), & \text{if } \lim_{\theta \to \infty} \frac{\alpha(\theta)}{\theta} = \infty. \end{cases}
$$

From Corollary 5.1, it follows that the selection cannot be detected at large deviation level for $\alpha(\theta) = o(\theta)$. In other words, a phase transition occurs at the critical scale $\theta$ which is different from the critical scale in (44).

In a recent paper, Joyce and Gao [7] studied the infinite-alleles model with homozygote advantage. This corresponds to choosing $H(\mathbf{p}) = \phi_2(\mathbf{p})$ in (39). A critical phenomenon is shown to exist in this case. They even obtained the following corollary to Theorem 5.2.

COROLLARY 5.2. *Choosing $H(\mathbf{p}) = \phi_2(\mathbf{p})$ in (39), then the family $\Xi_{\alpha,\theta}^H$ satisfies an LDP with speed $\theta$ and rate function*

$$
(46) \quad S_\alpha(\mathbf{p}) = \begin{cases} S(\mathbf{p}), & \text{if } \lim_{\theta \to \infty} \frac{\alpha(\theta)}{\theta} = 0, \\ -cH(\mathbf{p}) + S(\mathbf{p}), & \text{if } \lim_{\theta \to \infty} \frac{\alpha(\theta)}{\theta} = c \leq c_0, \\ \log\left(\frac{1 - \sqrt{1 - 2/c}}{2}\right) \\ \quad + c\left(\frac{1 + \sqrt{1 - 2/c}}{2}\right)^2 \\ \quad - cH(\mathbf{p}) + S(\mathbf{p}), & \text{if } \lim_{\theta \to \infty} \frac{\alpha(\theta)}{\theta} = c > c_0, \end{cases}
$$

*where $c_0 > 2$ solves the equation*

$$
\log\left(\frac{1 - \sqrt{1 - 2/c}}{2}\right) + c\left(\frac{1 + \sqrt{1 - 2/c}}{2}\right)^2 = 0.
$$

**Acknowledgments.** We thank the referees for their comments and suggestions, and for explaining the historical development of the study of the asymptotic behavior of the Poisson–Dirichlet distribution for large $\theta$. We also thank Paul Joyce for informing us of the result in Corollary 5.2.

SCHOOL OF MATHEMATICS AND STATISTICS  
CARLETON UNIVERSITY  
OTTAWA, ONTARIO  
CANADA K1S 5B6  
E-MAIL: ddawson@math.carleton.ca  

DEPARTMENT OF MATHEMATICS AND STATISTICS  
MCMASTER UNIVERSITY  
HAMILTON, ONTARIO  
CANADA L8S 4K1  
E-MAIL: shuifeng@univmail.cis.mcmaster.ca